\documentstyle[12pt,twoside]{article}
\newtheorem{theorem}{\bf Theorem}[section]
\newtheorem{proposition}[theorem]{\bf Proposition}

\newtheorem{corollary}[theorem]{\bf Corollary}

\input{amssym}

\date{}
\begin{document}
\title{{\Large\bf On $(p, q)-$centralizers of certain Banach algebras}}

\author{{\normalsize\sc M. J. Mehdipour\footnote{Corresponding author}~ and N. Salkhordeh}}
\maketitle

{\footnotesize {\bf Abstract.}   Let $A$ be an algebra with a right identity. In this paper, we study $(p, q)-$centralizers of $A$ and show that every $(p, q)-$centralizer of $A$ is a two-sided centralizer. In the case where, $A$ is normed algebra, we also prove that $(p, q)-$centralizers of $A$ are bounded.  Then, we apply the results for some group algebras and verify that $L^1(G)$ has a nonzero weakly compact $(p, q)-$centralizer if and only if $G$ is compact and the center of $L^1(G)$ is non-zero. Finally, we investigate $(p, q)-$Jordan centralizers of $A$ and determine them.
}
{}

\section{\normalsize\bf Introduction}

Let $T$ be a map on an algebra $ A $, and let $ p $ and $ q $ be distinct non-negative integers. Then $T$ is called a $\textit{$(p, q)-$centralizer}$ if for every $a, b \in A $
\begin{eqnarray*}\label{t1}
(p+q) T(ab) = p T(a)b + q a T(b)
\end{eqnarray*}
for all $ a \in A $. We denote by $ C_{p,q}(A) $ the space of all $(p, q)-$centralizers of $ A $. Every $(1,0)-$centralizer is called a \emph{left  centralizer} and every $(0,1)-$centralizer is called a \emph{right  centralizer}. Moreover,  $T$ is said to be a $\textit{two-sided centralizer }$ if $T$ is both a left centralizer and a right centralizer. Let $C_{ts}(A)$ be the set of all two-sided centralizers of $A$. It is easy to see that if $a\in A$, then the linear mapping $\rho_a: A\rightarrow A$ defined by $\rho_a(b)=ba$ is a right centralizer. Also, the linear mapping $\lambda_a: A\rightarrow A$ defined by $\lambda_a(b)=ab$ is a left centralizer. Note that if $a\in Z(A)$, the center of $A$, then $\rho_a=\lambda_a\in C_{ts}(A)$.
An additive mapping $T$ on an algebra $A$ is a $ \textit{$(p, q)-$Jordan centralizer} $ if
\begin{eqnarray*}\label{t2}
(p+q) T(a^{2}) = p T(a)a + q a T(a)
\end{eqnarray*}
The set of all $(p, q)-$Jordan centralizers of $A$ is denoted by $C_{p, q}^J(A)$.Clearly,
$$
C_{ts}(A)\subseteq C_{p, q}(A)\subseteq C_{p, q}^J(A).
$$
Vukman [19] introduced the notion of $(p, q)-$Jordan centralizers. For a prime ring $R$ with a suitable characteristic and positive integers $p$ and $q$, he proved that
$$
C_{p, q}^J(R)= C_{ts}(R),
$$
when $Z(R)$ is non-zero. Kosi-Ulbl and Vukman [8] proved this result for semiprime rings; see [4, 6] for generalized $(p, q)-$Jordan centralizers.

In this paper, we investigate $(p, q)-$centralizers and $(p, q)-$Jordan centralizers of algebras. In Section 2, we show that $(p, q)-$centralizers of an algebra $A$ with a right identity are two-sided centralizers. In the case where, $A$ is a normed algebra, every $(p, q)-$centralizer of $A$is bounded.  We prove that every bounded $(p, q)-$centralizer of an algebra with a bounded approximate identity is a two-sided centralizer. In Section 3, we give some results on the range of $(p, q)-$centralizers. In Section 4, we investigate  $(p, q)-$centralizers of some group algebras and determine the set of all their $(p, q)-$centralizers. For a locally compact group $G$, we show that $L^1(G)$ has a non-zero weakly compact $(p, q)-$centralizer if and only if $G$ is compact and the center of $L^1(G)$ is non-zero. In Section 5, we give a characterization of $(p, q)-$Jordan centralizers of an algebra with a right identity. We also show that every $(p, q)-$Jordan centralizer of a commutative algebra is always a two-sided centralizer.

{\bf Convention:} In this paper, $p$ and $q$ are distinct  positive integers.

\section{\normalsize\bf $(p, q)-$centralizers on algebras}

First, we characterize $ (p, q)-$centralizers of algebras.

\begin{theorem}\label{mj4} Let $A$ be an algebra with a right identity $ u$. If  $T: A\rightarrow A$ is a map, then the following assertions are equivalent.

\emph{(a)} $T\in C_{p, q}(A)$.

\emph{(b)} $ T=\rho_{T(u)}\in C_{ts}(A)$.

\emph{(c)} For every $a, b\in A$, $aT(b)=T(a)b$.\\
In this case, $T$ is linear.
\end{theorem}
{\it Proof.} Let $T\in C_{p, q}(A)$. Then for every $ a \in A $, we have
\begin{eqnarray}
(p+q) T(a) &=& (p+q)T(au)  \nonumber\\
&=& p T(a) + q a T(u). \nonumber
\end{eqnarray}
Thus
\begin{eqnarray}\label{r1}
T(a) = a T(u).
\end{eqnarray}
It follows that
\begin{eqnarray}
(p+q) abT(u) &=& (p+q)T(ab)  \nonumber\\
&=& p T(a)b + q a T(b)  \nonumber\\
&=& p aT(u)b + q ab T(u). \nonumber
\end{eqnarray}
Hence $ abT(u) = aT(u)b$. From this and (1) we see that $ aT(b) = T(a)b $.
Thus for every $ a,b \in A $,
\begin{eqnarray}
(p+q)T(ab) &=& p T(a)b + q a T(b)  \nonumber\\
&=& (p+q)T(a)b. \nonumber
\end{eqnarray}
So, $ T(ab)=T(a)b = a T(b). $ That is, $T\in C_{ts}(A)$. Therefore, (a)$\Rightarrow$(b). The implications (b)$\Rightarrow$(c)$\Rightarrow$(a) are clear.$\hfill\square$\\

For a normed algebra $A$, let $C_{p, q}^b(A)$ be the space of all bounded $(p, q)-$centralizers of $A$.

\begin{corollary}
Let $A$ be a normed algebra with a right identity. Then $C_{p, q}(A)$ is a Banach algebra and
$$C_{p, q}(A)=C_{p, q}^b(A)=\{\rho_a: a\in A\}.$$
\end{corollary}

Let $A$ and $B$ be two subsets of an algebra $\frak{A}$. We denote by $Z(A, B)$ the set of all $a\in A$ such that $ab=ba$ for all $b\in B$. Note that if $A=B$, then $Z(A, B)=Z(A)$.

\begin{corollary}\label{mj9}
Let $ A $ be an algebra with identity $ 1_{A} $. Then the following statements hold.

\emph{(i)} If $T\in C_{p, q}(A)$, then $ T(1_{A}) \in Z(A) $ and $ T= \rho_{T(1_{A})} = \lambda_{T(1_{A})}$.

\emph{(ii)} If $A$ is also a normed algebra, then $ C_{p,q}(A)\cong Z(A)$, where ``$\cong$" denotes isometrically isomorphic as Banach algebras.
\end{corollary}
{\it Proof.} (i) Let $T\in C_{p, q}(A)$. In view of Theorem 2.1, $T\in C_{ts}(A)$. Thus for every $ a \in A $, we have
$$ T(a) = T(a1_{A}) = a T(1_{A}) $$
and
$$ T(a) = T(1_{A}a) = T(1_{A} )a.$$
Hence $ a T(1_{A}) = T(1_{A} )a $ and so $ T(1_{A} ) \in Z(A)$.

(ii) Let $A$ be an algebra with the identity $1_A$. Then $ C_{p,q}(A) $ is a  subalgebra of algebra of all bounded linear operators on $A$. We define the function $\Gamma: C_{p,q}(A)\rightarrow Z(A)$ by $\Gamma(T)=T(1_A)$. It is clear that $\Gamma$ is linear. Let  $T\in C_{p,q}(A)$. Then for every $a\in A$,
$$
\|T(a)\|=\|aT(1_A)\|\leq\|a\|\|T(1_A)\|.
$$
This implies that $$\|T\|=\|T(1_A)\|=\|\Gamma(T)\|.$$ Therefore, $\Gamma$ is isometric.

If $T_1, T_2\in C_{p,q}(A)$, then
\begin{eqnarray*}
\Gamma(T_1)\Gamma(T_2)&=&T_1(1_A)T_2(1_A)\\
&=&T_1(1_AT_2(1_A))\\
&=&T_1T_2(1_A)\\
&=&\Gamma(T_1T_2).
\end{eqnarray*}
Hence $\Gamma$ is an algebra homomorphism.  Also, for every $c\in Z(A)$, the linear map $\rho_c$ is an element of $C_{p,q}(A)$. So $\Gamma$ is isomorphic.
$\hfill\square$\\

Let $A$ be a Banach algebra. It is well-known that the second conjugate of $A$, denoted by $A^{**}$, with the first Arens product ``$\diamond$" defined by
$$
\langle F\diamond H, f\rangle=\langle F, Hf\rangle
$$
is a Banach algebra, where $\langle Hf, a\rangle=\langle H, fa\rangle$, in which
$\langle fa, b\rangle=\langle f, ab\rangle$ for all $F, H\in A^{**}$, $f\in A^*$ and $a, b\in A$; see for example [5].

\begin{theorem}\label{mj6} Let $ A $ be a Banach algebra with a bounded approximate identity. Then $C_{p, q}^b(A)$ is a closed subspace of $C_{ts}(A)$.
\end{theorem}
{\it Proof.} Let  $T\in C^b_{p, q}(A)$. Then
\begin{eqnarray*}
(p+q) \langle T^{*}(f) a, x \rangle  &=& (p+q) \langle T^{*}(f), ax \rangle \\
 &=&  \langle (p+q) f , T(ax) \rangle \\
 &=&  \langle f , p T(a) x + q a T(x)  \rangle\\
 &=& p \langle f , T(a) x \rangle  + q \langle f, a T(x)  \rangle \\
 &=& \langle p f T(a), x \rangle  +  \langle q T^{*}(fa), x  \rangle
\end{eqnarray*}
for all $ a,x \in A  $. So, if $ H \in A^{**} $, then
\begin{eqnarray}
(p+q) \langle HT^{*}(f), a \rangle  &=& (p+q) \langle H,  T^{*}(f) a \rangle  \nonumber \\
 &=&  \langle H, p f d (a)+ q T^{*}(fa) \rangle  \nonumber \\
 &=& p \langle T^{*} (Hf), a \rangle + q \langle  T^{**} (H) f, a \rangle .  \nonumber
\end{eqnarray}
Thus for every $ F \in A^{**} $ we obtain
\begin{eqnarray}
\langle (p+q) T^{**} (F\diamond H), f \rangle &=& (p+q) \langle F\diamond H, T^{*}(f) \rangle  \nonumber \\
 &=& (p+q) \langle F, HT^{*}(f) \rangle  \nonumber \\
 &=&  \langle F, p T^{*} (Hf) + q T^{**} (H) f  \rangle  \nonumber \\
 &=&  \langle p T^{**} (F)H + q F T^{**} (H), f  \rangle .  \nonumber
\end{eqnarray}
That is, $ T^{**} $ is a $(p, q)-$centralizer. Since $A^{**}$ has a right identity, by Theorem 2.1,  $ T^{**}\in C_{ts}(A^{**}) $. Therefore, $T\in C_{ts}(A)$.
$\hfill\square$

\section{\normalsize\bf The range of $(p, q)-$centralizers}

\begin{theorem} Let $ A $ be a normed algebra with a right identity $u$, and let $T\in C_{p, q}(A)$. Then the following assertions are equivalent.

\emph{(a)} $T$ maps $A$ into $uA$.

\emph{(b)} $T=\lambda_{a_0}$ for some $a_0\in A$.

\emph{(c)} $T=\lambda_{T(u)}$.

\emph{(d)} $T(u)\in Z(A)$.
\end{theorem}
{\it Proof.} Assume that $T\in C_{p, q}(A)$ maps $A$ into $uA$. Then by Theorem 2.1, for every $a\in A$ we have
$$
uT(a-ua)=T(u)(a-ua)=0.
$$
Thus
$$
T(a-ua)\in\hbox{ran}(A)\cap uA=\{0\}.
$$
Hence $T(a)=T(ua)=T(u)a=\lambda_{T(u)}(a).$ Thus  (a)$\Rightarrow$(b).

Let $T=\lambda_{a_0}$ for some $a_0\in A$. Then for every $a\in A$, we have
\begin{eqnarray*}
a_0&=&T(u)=\rho_{T(u)}(u)\\
&=&uT(u)=T(u)u\\
&=&T(u).
\end{eqnarray*}
Hence (b)$\Rightarrow$(c). The implication (c)$\Rightarrow$(d) follows from Theorem 2.1. Finally, if $T(u)\in Z(A)$, then
$$
T(a)=aT(u)=T(u)a=uT(a)
$$
for all $a\in A$. Thus $T$ maps $A$ into $uA$. That is, (d)$\Rightarrow$(a).
$\hfill\square$

\begin{proposition}\label{n2}
Let  $A$ be an algebra and $T\in C_{p, q}(A)$. Then $T^2=0$ if and only if the range of $T$ is nilpotent. In this case, $T$ maps $A$ into the radical of $A$.
\end{proposition}
{\it Proof.} For every  $ a\in A $, we have
\begin{eqnarray*}\label{1212}
0 &=& (p+q)T^{2}(a^2)  \nonumber\\
&=& T(p T(a)a + q a T(a)) \\
&=& p^{2}  T^{2}(a) a + pq  T(a)T(a)+ pq  T(a) T(a)+ q^{2} a  T^{2}(a) \nonumber\\
&=& 2pq T(a)^2.     \nonumber
\end{eqnarray*}
 This proves the result.
$\hfill\square$\\

Let $\frak{A}$ be a Banach algebra with a bounded approximate identity such that every proper closed ideal of $\frak{A}$ is contained in a proper closed ideal with a bounded approximate identity. From Theorem 4.7 of [18] and Theorem 2.4 we have the following result.

\begin{proposition} Let $\frak{A}$ be as above. Then $T\in C_{p, q}(\frak{A})$ has a closed range if and only if there exist an idempotent $(p, q)-$centralizer $T_1$ of $\frak{A}$ and an invertible $(p, q)-$centralizer $T_2$ of $\frak{A}$ such that $T=T_1\circ T_2=T_2\circ T_1$.
\end{proposition}

\section{\normalsize\bf $(p, q)-$centralizers on group algebras}

Let $G$ be a locally compact group with a left Haar measure. Let $L^1(G)$ be the group algebra of $G$. Then $L^1(G)$ with the convolution product ``$\ast$" and the norm $\|.\|_1$ is a Banach algebra with a bounded approximate identity [5, 7]. Let $L^\infty(G)$ be the usual Lebesgue space as defined in [7] and $L_0^\infty(G)$ be the subspace of $L^\infty(G)$ consisting of all functions  that vanish at infinity. Then $L^\infty(G)^*$ and $L_0^\infty(G)^*$ are Banach  algebras with the first Arens product. One can prove that $L^\infty(G)^*$ and $L_0^\infty(G)^*$ have right identities [5, 9]; for more study see [1, 2, 11-14].  Let $M(G)$ be the measure algebra of $G$. Then $M(G)$ with the convolution product is a unital Banach algebra and $M(G)\cong C_0(G)^*$, where $C_0(G)$ is the space of all complex-valued continuous functions on $G$ that vanish at infinity [7].  Finally, let $C_b(G)$ be the space of all bounded continuous functions on $G$, and let  $LUC(G)$ be the space of all $f\in C_b(G)$ such that the mapping $x\mapsto  l_xf$ from $G$ into $C_b(G)$ is continuous, where $l_xf(y)=f(xy)$ for all $y\in G$. Let us remark that $LUC(G)^*$ with the product ``$\cdot$" defined by
$$
\langle m\cdot n, f\rangle=\langle m, nf\rangle\quad(m, n\in LUC(G)^*, f\in LUC(G)),
$$
where
$$
\langle nf, x\rangle=\langle n, l_xf\rangle\quad(x\in G)
$$
is a unital Banach algebra [5]. These facts  together with Corollary 2.3 show that $C_{p, q}(\frak{A})$ is a Banach algebra and $C_{p, q}(\frak{A})\cong Z(\frak{A})$, where ``$\cong$" denotes isometrically isomorphic as Banach algebras and $\frak{A}$ is one of the Banach algebras $M(G)$ or $LUC(G)^*$.

\begin{theorem}\label{mah} Let $G$ be a locally compact group.
Then the following statements hold.

\emph{(i)} $C_{p, q}^b(L^1(G))=\{\rho_\mu: \mu\in Z(M(G), L^1(G))\}.$

\emph{(ii)} If $T\in C_{p, q}(L_0^\infty(G)^*)$, then $T=\rho_\mu$ for some $\mu\in Z(M(G), L^1(G)).$
\end{theorem}
{\it Proof.} Let $\frak{B}$ be one of the Banach algebras $L^1(G)$ or $L_0^\infty(G)^*$.  Assume that $T\in C_{p, q}^b(\frak{B})$. Then by Theorems 2.1 and 2.4, $T\in C_{ts}(\frak{B})$. So $T$ is a right multiplier on $\frak{B}$. Note that if $\frak{B}=L_0^\infty(G)^*$, then $T=\rho_n$, where $n=T(u)$ and $u$ is a right identity for $L_0^\infty(G)^*$ with $\|u\|=1$. By Lemma 2.2 of [9],  $T=\rho_\mu$ for some $\mu\in M(G)$. If $\frak{B}=L^1(G)$, then by [20], $T=\rho_\mu$ for some $\mu\in M(G)$. Thus $T=\rho_\mu$ on $\frak{B}$. Now, let $(e_\alpha)_{\alpha\in\Lambda}$ be a bounded approximate identity of $L^1(G)$. Then for every $\phi\in L^1(G)$ and $\alpha\in\Lambda$, we have
\begin{eqnarray*}
e_\alpha\ast\phi\ast\mu&=&T(e_\alpha\ast\phi)\\
&=&T(e_\alpha)\ast\phi=e_\alpha\ast\mu\ast\phi.
\end{eqnarray*}
Since $L^1(G)$ is an ideal of $M(G)$, it follows that $\phi\ast\mu=\mu\ast\phi$. Thus $\mu\in Z(M(G), L^1(G))$.
$\hfill\square$\\

Let $C_{p, q}^w(A)$ be the space of all weakly compact $(p, q)-$centralizers of $A$.

\begin{corollary}\label{ab} Let $ G $ be a compact group. Then the following statements hold.

\emph{(i)} $C_{p, q}^w(L^1(G))=\{\rho_\phi: \phi\in Z(L^1(G))\}.$

\emph{(ii)} $C_{p, q}^w(M(G))=\{\rho_\phi: \phi\in Z(L^1(G), M(G))\}.$

\emph{(iii)} $C_{p, q}^w(L^\infty(G)^*)=\{\rho_\phi: \phi\in Z(L^1(G))\}$.

\emph{(iv)} $C_{p, q}^w(LUC(G)^*)=\{\rho_\phi: \phi\in Z(L^1(G))\}.$
\end{corollary}
{\it Proof.} (i) Let $T\in C_{p, q}^w(L^1(G))$. Then $T=\rho_\mu$ for some $\mu\in Z(M(G), L^1(G))$. Since $G$ is compact, every right centralizer of $L^1(G)$ is weakly compact and it is of the form $\rho_\phi$ for some $\phi\in L^1(G)$; see [3]. Thus $\mu\in L^1(G)$ and so $\mu\in Z(L^1(G))$.

(ii) Let $T\in C_{p, q}^w(M(G))$. Then $T=\rho_\mu$ for some $\mu\in Z(M(G))$. Thus $\rho_\mu$ is a weakly compact right multiplier on $L^1(G)$ which implies that $\mu\in L^1(G)$. So, $\mu\in Z(L^1(G), M(G))$.

(iii) Let $T\in C_{p, q}^w(L^\infty(G)^*)$. Then $T$ is a weakly compact right multiplier on $L^\infty(G)^*$. So $T=\rho_\phi$ for some $\phi\in L^1(G)$; see [15]. A similar argument to the proof of Theorem 4.1 shows that $\phi\in Z(L^1(G))$.

(iv) This follows from Theorem 2.1 and Corollary 2.3 of [16].
$\hfill\square$\\

Now, we prove the main result of this section.

\begin{theorem} Let $G$ be a locally compact group. Then the following assertions are equivalent.

\emph{(a)} $C_{p, q}^w(L_0^{\infty}(G)^{*})\neq\{0\}$.

\emph{(b)} $C_{p, q}^w(L^1(G))\neq\{0\}$.

\emph{(c)} $C_{p, q}^w(L^{\infty}(G)^{*})\neq\{0\}$.

\emph{(d)} $C_{p, q}^w(LUC(G)^{*})\neq\{0\}$.

\emph{(e)} $G$ is compact and $Z(L^1(G))\neq\{0\}$.
\end{theorem}
{\it Proof.}
Let $T\in C_{p, q}^w(L_0^{\infty}(G)^{*})$ be a non-zero. By Theorem 2.1, $T=\rho_n$ for some $n\in L_0^{\infty}(G)^{*}$. Since $L^1(G)$ is an ideal of $L_0^{\infty}(G)^{*}$, it follows that $T|_{L^1(G)}\in C_{p, q}^w(L^1(G))$. From weak$^*$ density of $L^1(G)$ into $L_0^{\infty}(G)^{*}$ we infer that $T|_{L^1(G)}$ is non-zero. So (a)$\Rightarrow$(b).

Let $T\in C_{p, q}^w(L^1(G))$ be  non-zero. According to Theorem 2.4, $T$ is a non-zero weakly compact  right centralizer of $L^1(G)$. Hence $G$ is compact; see [17]. By Corollary 4.2, $Z(L^1(G))\neq\{0\}$. Thus (b)$\Rightarrow$(e).

Assume that $Z(L^1(G))\neq\{0\}$ and $G$ is compact. It follows from Corollary 4.2 (iii) that $C_{p, q}^w(L^{\infty}(G)^{*})\neq\{0\}$. Hence
(e)$\Rightarrow$(a).

By [10] and Corollary 4.2, the implication (c)$\Rightarrow$(e) holds. The converse follows from the fact that $L^{\infty}(G)^{*}=L_0^{\infty}(G)^{*}$ when $G$ is compact.

Finally, Theorem 2.1 of [16] together with Corollary 4.2 (iv) proves  the implication (d)$\Rightarrow$(e); also the implication (e)$\Rightarrow$(d) holds by Corollary 2.3 of [16] and Corollary 4.2 (iv).$\hfill\square$\\

Let us remark that if $G$ is a locally compact abelian group, then $Z(L^1(G))=L^1(G)$.

\begin{corollary} Let $G$ be a locally compact abelian group. Then the following assertions are equivalent.

\emph{(a)} $C_{p, q}^w(L_0^{\infty}(G)^{*})\neq\{0\}$.

\emph{(b)} $C_{p, q}^w(L^1(G))\neq\{0\}$.

\emph{(c)} $C_{p, q}^w(L^{\infty}(G)^{*})\neq\{0\}$.

\emph{(d)} $C_{p, q}^w(LUC(G)^{*})\neq\{0\}$.

\emph{(e)} $G$ is compact.
\end{corollary}

\section{\normalsize\bf $(p, q)-$Jordan centralizers}

We commence this section with the following result.

\begin{theorem}\label{mj2}
Let  $ A $ be a commutative algebra. Then $ C_{p, q}^J(A)=C_{1, 1}(A)=C_{ts}(A)$.
\end{theorem}
{\it Proof.} First, assume that $m, n$ are distinct positive integers and $T\in C_{m, n}(A)$. Then for every $ a,b \in A $, we have
\begin{eqnarray}
m T(a) b + n a T(b) &=& (m+n)T(ab)  \nonumber\\
&=& (m+n)T(ba)  \nonumber\\
&=& m T(b)a + n b T(a) \nonumber
\end{eqnarray}
So $ T(a) b = a T(b) $. Hence for every $ a,b \in A $,
$$ T(ab) = T(a)b = a T(b).$$
Thus $C_{m, n}(A)=C_{ts}(A)$ for all  distinct positive integers $m, n$. This shows that $$C_{1, 1}(A)=C_{ts}(A).$$
Now, let $T\in C_{p, q}^J(A)$. Then for every $ a,b \in A $
$$ (p+q) T(ab+ba) = pT(a)b+pT(b)a+qaT(b)+qbT(a).$$
Since $ A $ is commutative, it follows that
$$ 2T(ab) = T(a)b + a T(b).$$
Consequently, $T\in C_{1, 1}(A)$ and therefore, $ C_{p, q}^J(A)=C_{1, 1}(A)$.
$\hfill\square$\\

We now give a characterization of $(p, q)-$Jordan centralizers of an algebra with a right identity.

\begin{theorem}\label{mj11}
Let $ A $ be an algebra with a right identity $ u $. If $T\in C_{p, q}^J(A)$, then
$$ T(a)=(a-ua)T(u)+uT(a)$$
for all $ a \in A $.
\end{theorem}
{\it Proof.}
Let $T\in C_{p, q}^J(A)$. Then for every $ a,b \in A $
\begin{eqnarray}\label{t3}
(p+q) T(ab+ba) = pT(a)b+pT(b)a+qaT(b)+qbT(a).
\end{eqnarray}
Put $ a=b=u $ in (2). Then $ T(u) = uT(u) $. Taking $ b=u $ in (2), we get
\begin{eqnarray}\label{t4}
qT(a)+(p+q) T(ua) = pT(u)a+qaT(u)+quT(a).
\end{eqnarray}
If we set $ a=ua $ in (3), then
\begin{eqnarray}\label{t5}
(p+2q) T(ua) = pT(u)a+quaT(u)+quT(ua).
\end{eqnarray}
From (3) we infer that
\begin{eqnarray*}
quT(a)+(p+q) uT(ua) = pT(u)a+quaT(u)+quT(a).
\end{eqnarray*}
Thus
\begin{eqnarray*}\label{t6}
(p+q) uT(ua) = pT(u)a+quaT(u)
\end{eqnarray*}
for all $ a \in A $. This together with (4) shows that
\begin{eqnarray*}
(p+2q) T(ua) &=&  pT(u)a+quaT(u)+ \frac{q}{p+q}(quaT(u)+pT(u)a) \\
&=&  \frac{p(p+2q)}{p+q}T(u)a+ \frac{q(p+2q)}{p+q}uaT(u).
\end{eqnarray*}
Hence
\begin{eqnarray*}\label{t7}
(p+q) T(ua) = pT(u)a+quaT(u).
\end{eqnarray*}
From this and (3) we see that
$$ qT(a)+quaT(u)+pT(u)a=qaT(u)+pT(u)a+quT(a).$$
This implies that
$$ T(a)=(a-ua)T(u)+uT(a)$$
for all $ a \in A $.
$\hfill\square$\\

\begin{theorem}\label{mj11}
Let $ A $ be an algebra with identity $ 1_{A} $ and $T\in C_{p, q}^J(A)$. If $ T(1_{A}) \in Z(A) $, then
$T\in C_{ts}(A)$.
\end{theorem}
{\it Proof.}
Let $T\in C_{p, q}^J(A)$. Then for every $ a,b \in A $
\begin{eqnarray}\label{t030}
(p+q) T(ab+ba) = pT(a)b+pT(b)a+qaT(b)+qbT(a).
\end{eqnarray}
Put $ b=1_{A} $ in (5). Then
\begin{eqnarray*}\label{t8}
(p+q) T(a) = pT(1_{A})a+qaT(1_{A}).
\end{eqnarray*}
If $ T(1_{A}) \in Z(A) $, then $ T(a) = aT(1_{A})= T(1_{A})a  $. Thus
$$ T(ab) = abT(1_{A})=a T(b) = aT(1_{A})b = T(a)b$$
That is, $T\in C_{ts}(A)$.
$\hfill\square$

\vspace{2mm}

{\footnotesize
\noindent {\bf Mohammad Javad Mehdipour}\\
Department of Mathematics,\\ Shiraz University of Technology,\\
Shiraz
71555-313, Iran\\ e-mail: mehdipour@sutech.ac.ir\\
\noindent {\bf Narjes Salkhordeh}\\
Department of Mathematics,\\ Shiraz University of Technology,\\ Shiraz
71555-313, Iran,\\ n.salkhordeh@sutech.ac.ir


\footnotesize

\begin{thebibliography}{99}

\bibitem{am} M. H. Ahmadi Gandomani and M. J. Mehdipour, Generalized derivations on some convolution algebras. Aequationes Math. 92 (2018), no. 2, 223--241.

\bibitem{am1} M. H. Ahmadi Gandomani and M. J. Mehdipour, Jordan, Jordan right and Jordan left derivations on convolution algebras. Bull. Iranian Math. Soc. 45 (2019), no. 1, 189--204.

\bibitem{a} C. A. Akemann, Some mapping properties of the group algebras of a compact group. Pacific J. Math. 22 (1967), 1--8.

\bibitem{bdf} D. Bennis, B. Dhara and B. Fahid, More on the generalized $(m,n)$-Jordan derivations and centralizers on certain semiprime rings. Bull. Iranian Math. Soc. 47 (2021), no. 1, 217--224.

\bibitem{d} H. G. Dales, Banach algebras and Automatic Continuity, Clarendon Press, Oxford, 2000.

\bibitem{f} A. Fo\v{s}ner, Ajda. A note on generalized $(m,n)$-Jordan centralizers. Demonstratio Math. 46 (2013), no. 2, 257--262.

\bibitem{hr} E. Hewitt and K. Ross, Abstract Harmonic Analysis I, Springer-Verlag, New York, 1970.

\bibitem{kv} I. Kosi-Ulbl and J. Vukman, On $(m,n)$-Jordan centralizers of semiprime rings. Publ. Math. Debrecen 89 (2016), no. 1-2, 223--231.

\bibitem{lp} A. T. Lau and J. Pym, Concerning the second dual of the group algebra of a  locally compact group. J. London Math. Soc. (2) 41 (1990), no. 3, 445--460.

\bibitem{l} V. Losert, Weakly compact multipliers on group algebras. J. Funct. Anal. 213 (2004), no. 2, 466--472.

\bibitem{mmn} S. Maghsoudi, M. J. Mehdipour and R. Nasr-Isfahani, Compact right multipliers on a Banach algebra related to locally compact semigroups. Semigroup Forum 83 (2011), no. 2, 205--213.

\bibitem{mmo} M. J. Mehdipour and Gh. R. Moghimi, The existence of nonzero compact right multipliers and Arens regularity of weighted Banach algebras. Rocky Mountain J. Math. 52 (2022), no. 6, 2101--2112.

\bibitem{mn1} M. J. Mehdipour and R. Nasr-Isfahani, Completely continuous elements of Banach algebras related to locally compact groups. Bull. Austral. Math. Soc. 76 (2007), no. 1, 49--54.

\bibitem{mn2}  M. J. Mehdipour and R. Nasr-Isfahani, Compact left multipliers on Banach algebras related to locally compact groups. Bull. Aust. Math. Soc. 79 (2009), no. 2, 227--238.

\bibitem{mn3} M. J. Mehdipour and R. Nasr-Isfahani, Weakly compact multipliers on Banach algebras related to a locally compact group. Acta Math. Hungar. 127 (2010), no. 3, 195--206.

\bibitem{m} M. J. Mehdipour, Weakly completely continuous elements of the Banach algebra LUC$(G)^*$. J. Math. Ext. 8 (2014), no. 1, 1--10.

\bibitem{s} S. Sakai, Weakly compact operators on operator algebras. Pacific J. Math. 14 (1964), 659--664.

\bibitem{u} A. \"{U}lger, When is the range of a multiplier on a Banach algebra closed? Math. Z. 254 (2006), no. 4, 715--728.

\bibitem{v} J. Vukman, On $(m,n)$-Jordan centralizers in rings and algebras. Glas. Mat. Ser. III 45(65) (2010), no. 1, 43--53.

\bibitem{w} J. G. Wendel, Left centralizers and isomorphisms of group algebras. Pacific J. Math. 2 (1952), 251--261.

\end{thebibliography}
\end{document}